\documentclass[12pt]{article}
\usepackage{amsmath,amsfonts,amssymb,amsthm}
\usepackage{color}
\usepackage{colortbl}
\usepackage{array}
\usepackage{mathrsfs}
\usepackage{dcolumn}
\usepackage{longtable}
\usepackage{hhline}
\usepackage{bm}

\newcommand{\braced}[2]{\genfrac{\{}{\}}{0pt}{0}{#1}{#2}}

\theoremstyle{plain}
\newtheorem{thm}{Theorem}[section]

\newtheorem{cor}[thm]{Corollary}

\theoremstyle{definition}
\newtheorem{defn}[thm]{Definition}

\title{\bf More Properties on Multi Poly-Euler Polynomials}
\author{
{\large Hassan Jolany and Roberto B. Corcino}
}

\date{}

\begin{document}

\maketitle

\begin{abstract}
In this paper, we establish more properties of generalized poly-Euler polynomials with three parameters and we investigate a kind of symmetrized generalization of poly-Euler polynomials. Moreover, we introduce a more general form of multi poly-Euler polynomials and obtain some identities parallel to those of the generalized poly-Euler polynomials.   


\end{abstract}

\section{Introduction}
The nature of introducing Bernoulli numbers is parallel to that of Euler numbers. Bernoulli numbers have been introduced by Jacob Bernoulli (1655-1705) in his effort to describe the coefficients of the polynomial representation of the sum
$$S_m(n)=1^m+2^m+\ldots +n^m,$$
while Euler numbers have been introduced by Leonard Euler (1707-1783) in his desire to evaluate the alternating sum
$$A_n(m)=m^n-(m-1)^n+\ldots +(-1)^{m-1}1^n.$$
Moreover, Bernoulli and Euler numbers have been usually defined by means of the following generating functions with parallel structures 
\begin{align}
\frac{x}{e^x-1}&=\sum_{n=0}^{\infty}B_n\frac{x^n}{n!}\\
\frac{2}{e^x+1}&=\sum_{n=0}^{\infty}E_n\frac{x^n}{n!},
\end{align}
respectively. It is worth mentioning that Euler worked intensively on Bernoulli numbers and gave great contributions in the development of the numbers. This made him known as ``godfather" of Bernoulli numbers (see \cite{Benyi}).

\smallskip
In 1997, Kaneko \cite{Kaneko} introduced the poly-Bernoulli numbers $B^{(k)}_n$ by means of the following exponential generating function
$$\frac{{\rm Li}_k(1-e^{-x})}{1-e^{-x}}=\sum_{n=0}^{\infty}B^{(k)}_n\frac{x^n}{n!}$$  
where
$${\rm Li}_k(z)=\sum_{n=0}^{\infty}\frac{z^n}{n^k}.$$
Ohno and Sasaki \cite{Sasaki}, on the other hand, defined poly-Euler numbers as
$$\frac{{\rm Li}_k(1-e^{-4t})}{4t\cosh t}=\sum_{n=0}^{\infty}E^{(k)}_n\frac{t^n}{n!}$$
which have been recently extended by H. Jolany et al. \cite{Jolany3} in polynomial form as
\begin{equation}\label{polyform}
\frac{2{\rm Li}_k(1-e^{-t})}{1+e^{t}}e^{xt}=\sum_{n=0}^{\infty}E^{(k)}_n(x)\frac{t^n}{n!}.
\end{equation}
Further generalization and other properties of poly-Bernoulli and poly-Euler numbers and polynomials including their relations with other special numbers and functions  are found in \cite{Araci, Bayad, Bayad2, Brewbaker, Candel, Coppo, Hamahata, Jang, Kim}. 

\smallskip
In \cite{Jolany3}, the generalized poly-Euler polynomials with parameters $a$, $b$ and $c$ are defined by
\begin{equation}\label{genpolyeuler1}
\frac{2{\rm Li}_k(1-(ab)^{-t})}{a^{-t}+b^{t}}c^{xt}=\sum_{n=0}^{\infty}{E}^{(k)}_n(x;a,b,c)\frac{t^n}{n!}.
\end{equation}
Note that the poly-Euler polynomials in \cite{Bayad, Sasaki} can be deduced from (\ref{genpolyeuler1}) by replacing $t$ with $4t$ and taking $x=1/2$. Moreover, when $x=0$, (\ref{genpolyeuler1}) gives 
$${E}^{(k)}_n(0;a,b,c)={E}^{(k)}_n(a,b)$$
where
$$\frac{2{\rm Li}_k(1-(ab)^{-t})}{a^{-t}+b^{t}}=\sum_{n=0}^{\infty}{E}^{(k)}_n(a,b)\frac{t^n}{n!},$$
and when $a=1$ and $b=c=e$ with ${E}^{(k)}_n(x;1,e,e)={E}^{(k)}_n(x)$, we obtain equation (\ref{polyform}). However, only one identity has been obtained in \cite{Jolany3} for ${E}^{(k)}_n(x;a,b,c)$ which is given by
\begin{equation}\label{explicit1}
{E}^{(k)}_n(x;a,b,c)=\sum_{m=0}^n\sum_{j=0}^m\sum_{i=0}^j\frac{2(-1)^{m-j+i}}{j^k}\binom{j}{i}(x\ln c-(m-j+i+1)\ln a-(m-j+i+1)\ln b)^n.
\end{equation}

\smallskip
In the same paper by H. Jolany et. al \cite{Jolany3}, they defined certain multi poly-Euler polynomials as follows
\begin{equation}\label{genpolyeuler}
\frac{2{\rm Li}_{(k_1, k_2,\ldots, k_r)}(1-(ab)^{-t})}{(a^{-t}+b^{t})^r}e^{rxt}=\sum_{n=0}^{\infty}{E}^{(k_1, k_2,\ldots, k_r)}_n(x;a,b)\frac{t^n}{n!}.
\end{equation}
where
$${\rm Li}_{(k_1,k_2,\ldots, k_r)}(z)=\sum_{ 0\le m_1\le m_2\le\ldots\le m_r }\frac{z^{m_r}}{m_1^{k_1} m_2^{k_2}\ldots m_r^{k_r}}$$
is the generalization of poly-logarithm, also known as multiple zeta values, which have been given much attention recently but until now, their precise structure remains a mystery. Note that, when $r=1$, (\ref{genpolyeuler}) immediately yields (\ref{genpolyeuler1}). Several identities on ${E}^{(k_1, k_2,\ldots, k_r)}_n(x;a,b)$ have been obtained in \cite{Jolany3} including the recurrence relations and certain explicit formula. However, this explicit formula is limited only to the case where $a=1$ and $b=e$. In this present paper, more identities for ${E}^{(k)}_n(x;a,b,c)$ will be established and further generalization of ${E}^{(k_1, k_2,\ldots, k_r)}_n(x;a,b)$ will be investigated.

\section{Some Results on Generalized Poly-Euler\\ Polynomials}
The main objective of this section is to establish more identities for ${E}^{(k)}_n(x;a,b,c)$. First, let us consider an expression for ${E}^{(k)}_n(x;a,b,c)$ in terms of ${E}^{(k)}_i(a,b)$, $i=0, 1, \ldots, n$.

\bigskip
\begin{thm}\label{thm1}
The generalized poly Euler polynomials satisfy the following relation
\begin{equation}\label{eqnnew1}
{E}^{(k)}_n(x;a,b,c)=\sum_{i=0}^n\binom{n}{i}(\ln c)^{n-i}{E}^{(k)}_i(a,b)x^{n-i}
\end{equation}
\begin{proof}
Using (\ref{genpolyeuler1}), we have
\begin{eqnarray*}
\sum_{n=0}^{\infty}{E}^{(k)}_n(x;a,b,c)\frac{t^n}{n!}&=&\frac{2{\rm Li}_k(1-(ab)^{-t})}{a^{-t}+b^{t}}c^{xt}=e^{xt\ln c}\sum_{n=0}^{\infty}{E}^{(k)}_n(a,b)\frac{t^n}{n!}\\
&=&\sum_{n=0}^{\infty}\sum_{i=0}^n\frac{(xt\ln c)^{n-i}}{(n-i)!}{E}^{(k)}_i(a,b)\frac{t^{i}}{i!}\\
&=&\sum_{n=0}^{\infty}\left(\sum_{i=0}^n\binom{n}{i}(\ln c)^{n-i}{E}^{(k)}_i(a,b)x^{n-i}\right)\frac{t^{n}}{n!}.
\end{eqnarray*}
Comparing the coefficients of $\frac{t^{n}}{n!}$, we obtain the desired result.
\end{proof}
\end{thm}

The next identity gives a relation between ${E}^{(k)}_n(x;a,b,c)$ and ${E}^{(k)}_n(x)$.

\bigskip
\begin{thm}\label{thm2}
The generalized poly Euler polynomials satisfy the following relation
\begin{equation}\label{eqnnew2}
{E}^{(k)}_n(x;a,b,c)=(\ln a+\ln b)^n{E}^{(k)}_n\left(\frac{x\ln c+\ln a}{\ln a+\ln b}\right)
\end{equation}
\begin{proof}
Using (\ref{genpolyeuler1}), we have
\begin{eqnarray*}
\sum_{n=0}^{\infty}{E}^{(k)}_n(x;a,b,c)\frac{t^n}{n!}&=&\frac{2{\rm Li}_k(1-(ab)^{-t})}{a^{-t}(1+(ab)^{t})}e^{xt\ln c}\\
&=&2e^{\frac{x\ln c+\ln a}{\ln ab}t\ln ab}\frac{{\rm Li}_k(1-e^{-t\ln ab})}{1+e^{t\ln ab}}\\
&=&\sum_{n=0}^{\infty}(\ln a+\ln b)^n{E}^{(k)}_n\left(\frac{x\ln c+\ln a}{\ln a+\ln b}\right)\frac{t^{n}}{n!}.
\end{eqnarray*}
Comparing the coefficients of $\frac{t^{n}}{n!}$, we obtain the desired result.
\end{proof}
\end{thm}

\bigskip
\begin{thm}\label{thm3}
%The generalized poly-Euler polynomials are Appell polynomials in the sense that
The generalized poly-Euler polynomials satisfy the following relation
\begin{equation}\label{eqnnew3}
\frac{d}{dx}{E}^{(k)}_{n+1}(x;a,b,c)=(n+1)(\ln c){E}^{(k)}_{n}(x;a,b,c)
\end{equation}
\begin{proof}
Using (\ref{genpolyeuler1}), we have
\begin{eqnarray*}
\sum_{n=0}^{\infty}\frac{d}{dx}{E}^{(k)}_{n}(x;a,b,c)\frac{t^{n}}{n!}&=&\frac{2t(\ln c){\rm Li}_k(1-(ab)^{-t})}{(a^{-t}+b^{t})}e^{xt\ln c}\\
\sum_{n=0}^{\infty}\frac{d}{dx}{E}^{(k)}_{n}(x;a,b,c)\frac{t^{n-1}}{n!}&=&\sum_{n=0}^{\infty}(\ln c){E}^{(k)}_{n}(x;a,b,c)\frac{t^{n}}{n!}.
\end{eqnarray*}
Hence, 
\begin{eqnarray*}
\sum_{n=0}^{\infty}\frac{1}{n+1}\frac{d}{dx}{E}^{(k)}_{n+1}(x;a,b,c)\frac{t^{n}}{n!}=\sum_{n=0}^{\infty}(\ln c){E}^{(k)}_{n}(x;a,b,c)\frac{t^{n}}{n!}.
\end{eqnarray*}
Comparing the coefficients of $\frac{t^{n}}{n!}$, we obtain the desired result.
\end{proof}
\end{thm}

The following corollary immediately follows from Theorem \ref{thm3} by taking $c=e$. For brevity, let us denote ${E}^{(k)}_{n}(x;a,b,e)$ by ${E}^{(k)}_{n}(x;a,b)$.

\bigskip
\begin{cor}\label{cor1}
The generalized poly-Euler polynomials are Appell polynomials in the sense that
\begin{equation}\label{eqnnew4}
\frac{d}{dx}{E}^{(k)}_{n+1}(x;a,b)=(n+1){E}^{(k)}_{n}(x;a,b)
\end{equation}
\end{cor}

\smallskip
Consequently, using the characterization of Appell polynomials \cite{Lee, Shohat, Toscano}, the following addition formula can easily be obtained.  

\bigskip
\begin{cor}\label{cor2}
The generalized poly-Euler polynomials satisfy the following addition formula
\begin{equation}\label{eqnnew5}
{E}^{(k)}_{n}(x+y;a,b)=\sum_{i=0}^n\binom{n}{i}{E}^{(k)}_i(x;a,b)y^{n-i}
\end{equation}
\end{cor}

\smallskip
\noindent Taking $x=0$ in formula (\ref{eqnnew5}) and using the fact that ${E}^{(k)}_{n}(0;a,b)={E}^{(k)}_{n}(a,b)$, Corollary \ref{cor2} gives formula (\ref{eqnnew1}) in Theorem \ref{thm1} with $c=e$.

\smallskip
Furthermore, using the fact that ${E}^{(k)}_{n}(x;a,b)$ satisfies the relation
\begin{equation}\label{genpolyeuler11}
\sum_{n=0}^{\infty}{E}^{(k)}_n(x;a,b)\frac{t^n}{n!}=\frac{2{\rm Li}_k(1-(ab)^{-t})}{a^{-t}+b^{t}}e^{xt},
\end{equation}
one can easily derive the following identities by manipulating the right-hand side of (\ref{genpolyeuler11}) and by appropriate application of the Cauchy's rule for the product of power series
\begin{align*}
{E}^{(k)}_{n}(x;a,b)&=\sum_{m=0}^{\infty}\frac{1}{m!}\sum_{l=m}^n\braced{l}{m}\binom{n}{l}{E}^{(k)}_{n-l}(-m;a,b)(x)^{(m)}\\
&=\sum_{m=0}^{\infty}\sum_{l=m}^n\braced{l}{m}\binom{n}{l}{E}^{(k)}_{n-l}(0;a,b)(x)_{m}\\
&=\sum_{m=0}^{\infty}\binom{n}{m}\sum_{l=0}^{n-m}\frac{\binom{n-m}{l}}{\binom{l+s}{l}}\braced{l+s}{s}{E}^{(k)}_{n-m-l}(0;a,b)B^{(s)}_m(x)\\
&=\sum_{m=0}^{\infty}\frac{\binom{n}{m}}{(1-\lambda)^s}\sum_{j=0}^s\binom{s}{j}(-\lambda)^{s-j}{E}^{(k)}_{n-m}(j;a,b)H^{(s)}_m(x;\lambda),
\end{align*}
where $(x)^{(n)}=x(x+1)\ldots (x+n-1)$, $(x)_n=x(x-1)\ldots (x-n+1)$,
$$\left(\frac{t}{e^t-1}\right)^se^{xt}=\sum_{n=0}^{\infty}B^{(s)}_n(x)\frac{t^n}{n!}\mbox{  and  }\left(\frac{1-\lambda}{e^t-\lambda}\right)^se^{xt}=\sum_{n=0}^{\infty}H^{(s)}_n(x;\lambda)\frac{t^n}{n!}.$$

\smallskip
Now, let us consider the following definition which contains certain symmetrized generalization of poly-Euler polynomials with parameters $a$, $b$ and $c$.

\begin{defn}\label{defn4}
For $m, n\ge0$, we define
\begin{equation}\label{symgen}
D^{(m)}_n(x,y;a,b,c)=\frac{1}{(\ln a+\ln b)^n}\sum_{k=0}^m\binom{m}{k}{E}^{(-k)}_n(x;a,b,c)\left(\frac{y\ln c+\ln a}{\ln a +\ln b}\right)^{m-k}.
\end{equation}
\end{defn}

The following theorem contains the double generating function for $D^{(m)}_n(x,y;a,b,c)$.

\bigskip
\begin{thm}\label{thm4}
For $n,m\ge0$, we have
\begin{equation}\label{eqnnew6}
\sum_{n=0}^{\infty}\sum_{m=0}^{\infty}D^{(m)}_n(x,y;a,b,c)\frac{t^n}{n!}\frac{u^m}{m!}=\frac{2e^{\left(\frac{y\ln c+\ln a}{\ln a +\ln b}\right)u}e^{\left(\frac{x\ln c+\ln a}{\ln a +\ln b}\right)t}e^{t+u}\left(1-e^{-t}\right)}{(e^t+1)(e^t+e^u-e^{t+u})}.
\end{equation}
\begin{proof}
$$\sum_{n=0}^{\infty}\sum_{m=0}^{\infty}D^{(m)}_n(x,y;a,b,c)\frac{t^n}{n!}\frac{u^m}{m!}\qquad\qquad\qquad\qquad\qquad\qquad\qquad\qquad\qquad\qquad\qquad\qquad\qquad\qquad\qquad\qquad$$
\begin{eqnarray*}
&=&\sum_{n=0}^{\infty}\sum_{m=0}^{\infty}\frac{1}{(\ln a+\ln b)^n}\sum_{k=0}^m{E}^{(-k)}_n(x;a,b,c)\left(\frac{y\ln c+\ln a}{\ln a +\ln b}\right)^{m-k}\frac{t^n}{n!}\frac{u^m}{k!(m-k)!}\\
&=&\sum_{n=0}^{\infty}\sum_{k=0}^{\infty}\sum_{m\ge k}\frac{1}{(\ln a+\ln b)^n}{E}^{(-k)}_n(x;a,b,c)\left(\frac{y\ln c+\ln a}{\ln a +\ln b}\right)^{m-k}\frac{t^n}{n!}\frac{u^m}{k!(m-k)!}.
\end{eqnarray*}
Replacing $m-k$ with $l$ and using Theorem \ref{thm2}, we obtain
$$\sum_{n=0}^{\infty}\sum_{m=0}^{\infty}D^{(m)}_n(x,y;a,b,c)\frac{t^n}{n!}\frac{u^m}{m!}\qquad\qquad\qquad\qquad\qquad\qquad\qquad\qquad\qquad\qquad\qquad\qquad\qquad\qquad\qquad\qquad$$
\begin{eqnarray*}
&=&\sum_{n=0}^{\infty}\sum_{k=0}^{\infty}\sum_{l=0}^{\infty}\frac{1}{(\ln a+\ln b)^n}{E}^{(-k)}_n(x;a,b,c)\left(\frac{y\ln c+\ln a}{\ln a +\ln b}\right)^{l}\frac{t^n}{n!}\frac{u^k}{k!}\frac{u^l}{l!}\\
&=&e^{\left(\frac{y\ln c+\ln a}{\ln a +\ln b}\right)u}\sum_{n=0}^{\infty}\sum_{k=0}^{\infty}{E}^{(-k)}_n\left(\frac{x\ln c+\ln a}{\ln a+\ln b}\right)\frac{t^n}{n!}\frac{u^k}{k!}\\
&=&e^{\left(\frac{y\ln c+\ln a}{\ln a +\ln b}\right)u}\sum_{k=0}^{\infty}\frac{2{\rm Li}_k(1-e^{-t})}{1+e^{-t}}e^{\left(\frac{x\ln c+\ln a}{\ln a+\ln b}\right)t}\frac{u^k}{k!}\\
&=&e^{\left(\frac{y\ln c+\ln a}{\ln a +\ln b}\right)u}e^{\left(\frac{x\ln c+\ln a}{\ln a +\ln b}\right)t}\sum_{k=0}^{\infty}\sum_{n=0}^{\infty}E^{(-k)}_n(0)\frac{t^n}{n!}\frac{u^k}{k!}.
\end{eqnarray*}
Note that
\begin{eqnarray*}
\sum_{k=0}^{\infty}\sum_{n=0}^{\infty}E^{(-k)}_n(0)\frac{t^n}{n!}\frac{u^k}{k!}&=&\sum_{k=0}^{\infty}\sum_{m>0}\frac{2\left(1-e^{-t}\right)^mm^k}{1+e^{t}}\frac{u^k}{k!}\\
&=&\sum_{m>0}\frac{2\left(1-e^{-t}\right)^m}{1+e^{t}}\sum_{k=0}^{\infty}\frac{(mu)^k}{k!}\\
&=&\frac{2}{1+e^{t}}\sum_{m>0}\left(1-e^{-t}\right)^me^{mu}\\
&=&\frac{2}{1+e^{t}}\frac{\left(1-e^{-t}\right)e^{u}}{1-\left(1-e^{-t}\right)e^{u}}=\frac{2e^{t+u}\left(1-e^{-t}\right)}{(e^t+1)(e^t+e^u-e^{t+u})}.
\end{eqnarray*}
Thus,
$$\sum_{n=0}^{\infty}\sum_{m=0}^{\infty}D^{(m)}_n(x,y;a,b,c)\frac{t^n}{n!}\frac{u^m}{m!}=\frac{2e^{\left(\frac{y\ln c+\ln a}{\ln a +\ln b}\right)u}e^{\left(\frac{x\ln c+\ln a}{\ln a +\ln b}\right)t}e^{t+u}\left(1-e^{-t}\right)}{(e^t+1)(e^t+e^u-e^{t+u})}.$$
\end{proof}
\end{thm}

The following is an explicit formula for $D^{(m)}_n(x,y;a,b,c)$.

\bigskip
\begin{thm}\label{thm5}
For $n,m\ge0$, we have
{\small
\begin{eqnarray*}
D^{(m)}_n(x,y;a,b,c)&=&2\sum_{j=0}^{\infty}(j!)^2\left(\sum_{l=0}^n\sum_{i=0}^{\infty}(-1)^i\frac{\left(\ln c^{x}a^{i+2}b^{i+1}\right)^{n-l}-\left(\ln c^{x}a^{i+1}b^{i}\right)^{n-l}}{(\ln a +\ln b)^{n-l}}\binom{n}{l}\braced{l}{j}\right)\times\\
&&\;\;\;\;\;\times\left(\sum_{r=0}^m\left(\frac{y\ln c+2\ln a+\ln b}{\ln a +\ln b}\right)^{m-r}\binom{m}{r}\braced{r}{j}\right)
\end{eqnarray*}
}
\begin{proof} Using Theorem \ref{thm4},
$$\sum_{n=0}^{\infty}\sum_{m=0}^{\infty}D^{(m)}_n(x,y;a,b,c)\frac{t^n}{n!}\frac{u^m}{m!}\qquad\qquad\qquad\qquad\qquad\qquad\qquad\qquad\qquad\qquad\qquad\qquad\qquad\qquad\qquad\qquad\qquad\qquad\qquad\qquad\qquad\qquad\qquad\qquad$$
\begin{eqnarray*}
&=&\frac{2e^{\left(\frac{y\ln c+\ln a}{\ln a +\ln b}\right)u}e^{\left(\frac{x\ln c+\ln a}{\ln a +\ln b}\right)t}e^{t+u}\left(1-e^{-t}\right)}{(e^t+1)(1- (e^t-1)(e^u-1))}\\
&=&2e^{\left(\frac{y\ln c+2\ln a+\ln b}{\ln a +\ln b}\right)u}e^{\left(\frac{x\ln c+2\ln a+\ln b }{\ln a +\ln b}\right)t}\left(1-e^{-t}\right)\sum_{i=0}^{\infty}(-1)^ie^{it}\sum_{j=0}^{\infty}(e^t-1)^j(e^u-1)^j\\
&=&2e^{\left(\frac{y\ln c+2\ln a+\ln b}{\ln a +\ln b}\right)u}\sum_{i=0}^{\infty}(-1)^ie^{\left(\frac{x\ln c+(i+2)\ln a+(i+1)\ln b }{\ln a +\ln b}\right)t}\sum_{j=0}^{\infty}(e^t-1)^j(e^u-1)^j\\
&&\;\;\;\;\;\;\;-2e^{\left(\frac{y\ln c+2\ln a+\ln b}{\ln a +\ln b}\right)u}\sum_{i=0}^{\infty}(-1)^ie^{\left(\frac{x\ln c+(i+1)\ln a+i\ln b }{\ln a +\ln b}\right)t}\sum_{j=0}^{\infty}(e^t-1)^j(e^u-1)^j.
\end{eqnarray*}
Applying the exponential generating function for Stirling numbers of the second kind \cite{Comtet}
$$\sum_{n=0}^{\infty}\braced{n}{k}\frac{t^n}{n!}=\frac{(e^t-1)^k}{k!},$$
we obtain
$$\sum_{n=0}^{\infty}\sum_{m=0}^{\infty}D^{(m)}_n(x,y;a,b,c)\frac{t^n}{n!}\frac{u^m}{m!}\qquad\qquad\qquad\qquad\qquad\qquad\qquad\qquad\qquad\qquad\qquad\qquad\qquad\qquad\qquad\qquad\qquad\qquad\qquad\qquad\qquad\qquad\qquad\qquad$$
\begin{eqnarray*}
&=&2\sum_{j=0}^{\infty}\left(j!\sum_{i=0}^{\infty}(-1)^i\sum_{n=0}^{\infty}\frac{\left(\frac{x\ln c+(i+2)\ln a+(i+1)\ln b }{\ln a +\ln b}\right)^nt^n}{n!}\sum_{m=0}^{\infty}\braced{m}{j}\frac{t^m}{m!}\right)\times\\
&&\;\;\;\;\;\times\left(j!\sum_{n=0}^{\infty}\frac{\left(\frac{y\ln c+2\ln a+\ln b}{\ln a +\ln b}\right)^nu^n}{n!}\sum_{m=0}^{\infty}\braced{m}{j}\frac{u^m}{m!}\right)\\
&&-2\sum_{j=0}^{\infty}\left(j!\sum_{i=0}^{\infty}(-1)^i\sum_{n=0}^{\infty}\frac{\left(\frac{x\ln c+(i+1)\ln a+i\ln b }{\ln a +\ln b}\right)^nt^n}{n!}\sum_{m=0}^{\infty}\braced{m}{j}\frac{t^m}{m!}\right)\times\\
&&\;\;\;\;\;\times\left(j!\sum_{n=0}^{\infty}\frac{\left(\frac{y\ln c+2\ln a+\ln b}{\ln a +\ln b}\right)^nu^n}{n!}\sum_{m=0}^{\infty}\braced{m}{j}\frac{u^m}{m!}\right)\\
&=&2\sum_{j=0}^{\infty}\left(j!\sum_{i=0}^{\infty}\sum_{l=0}^{\infty}\sum_{m=0}^l(-1)^i\left(\frac{x\ln c+(i+2)\ln a+(i+1)\ln b }{\ln a +\ln b}\right)^{l-m}\binom{l}{m}\braced{m}{j}\frac{t^l}{l!}\right)\times\\
&&\;\;\;\;\;\times\left(j!\sum_{p=0}^{\infty}\sum_{r=0}^p\left(\frac{y\ln c+2\ln a+\ln b}{\ln a +\ln b}\right)^{p-r}\binom{p}{r}\braced{r}{j}\frac{u^p}{p!}\right)\\
&&-2\sum_{j=0}^{\infty}\left(j!\sum_{i=0}^{\infty}\sum_{l=0}^{\infty}\sum_{m=0}^l(-1)^i\left(\frac{x\ln c+(i+1)\ln a+i\ln b }{\ln a +\ln b}\right)^{l-m}\binom{l}{m}\braced{m}{j}\frac{t^l}{l!}\right)\times\\
&&\;\;\;\;\;\times\left(j!\sum_{p=0}^{\infty}\sum_{r=0}^p\left(\frac{y\ln c+2\ln a+\ln b}{\ln a +\ln b}\right)^{p-r}\binom{p}{r}\braced{r}{j}\frac{u^p}{p!}\right)
\end{eqnarray*}
\begin{eqnarray*}
&=&\sum_{n=0}^{\infty}\sum_{m=0}^{\infty}\frac{t^n}{n!}\frac{u^m}{m!}2\sum_{j=0}^{\infty}(j!)^2\left(\sum_{l=0}^n\sum_{i=0}^{\infty}(-1)^i\left(\frac{x\ln c+(i+2)\ln a+(i+1)\ln b }{\ln a +\ln b}\right)^{n-l}\binom{n}{l}\braced{l}{j}\right)\times\\
&&\;\;\;\;\;\times\left(\sum_{r=0}^m\left(\frac{y\ln c+2\ln a+\ln b}{\ln a +\ln b}\right)^{m-r}\binom{m}{r}\braced{r}{j}\right)\\
&&-\sum_{n=0}^{\infty}\sum_{m=0}^{\infty}\frac{t^n}{n!}\frac{u^m}{m!}2\sum_{j=0}^{\infty}(j!)^2\left(\sum_{l=0}^n\sum_{i=0}^{\infty}(-1)^i\left(\frac{x\ln c+(i+1)\ln a+i\ln b }{\ln a +\ln b}\right)^{n-l}\binom{n}{l}\braced{l}{j}\right)\times\\
&&\;\;\;\;\;\times\left(\sum_{r=0}^m\left(\frac{y\ln c+2\ln a+\ln b}{\ln a +\ln b}\right)^{m-r}\binom{m}{r}\braced{r}{j}\right)
\end{eqnarray*}
Comparing the coefficients yields the desired result.
\end{proof}
\end{thm}

\section{Generalized Multi Poly-Euler Polynomials}
Let us define a more general form of multi poly-Euler polynomials.

\begin{defn}\label{multigenpolyeuler}
The generalized multi poly-Euler polynomials with parameters $a$, $b$ and $c$ are defined by
\begin{equation}\label{genpolyeuler2}
\frac{2{\rm Li}_{(k_1, k_2,\ldots, k_r)}(1-(ab)^{-t})}{(a^{-t}+b^{t})^r}c^{rxt}=\sum_{n=0}^{\infty}{E}^{(k_1, k_2,\ldots, k_r)}_n(x;a,b,c)\frac{t^n}{n!}.
\end{equation}
\end{defn}

\smallskip
\noindent In particular,
\begin{eqnarray*}
{E}^{(k_1, k_2,\ldots, k_r)}_n(x)&=&{E}^{(k_1, k_2,\ldots, k_r)}_n(x;1,e,e)\\
{E}^{(k_1, k_2,\ldots, k_r)}_n(a,b)&=&{E}^{(k_1, k_2,\ldots, k_r)}_n(0;a,b)
\end{eqnarray*}

\smallskip
The following theorem contains some identities for ${E}^{(k_1, k_2,\ldots, k_r)}_n(x;a,b,c)$ which can be derived using the process in deriving the identities in Theorems \ref{thm1} -- \ref{thm3}.

\bigskip
\begin{thm}\label{multithm1}
The generalized multi poly-Euler polynomials satisfy the following relations
\begin{eqnarray}
{E}^{(k_1, k_2, \ldots, k_r)}_n(x;a,b,c)&=&\sum_{i=0}^n\binom{n}{i}(r\ln c)^{n-i}{E}^{(k_1, k_2, \ldots, k_r)}_i(a,b)x^{n-i}\nonumber \\
{E}^{(k_1, k_2, \ldots, k_r)}_n(x;a,b,c)&=&(\ln a+\ln b)^n{E}^{(k_1, k_2, \ldots, k_r)}_n\left(\frac{rx\ln c+\ln a}{\ln a+\ln b}\right)\nonumber \\
\frac{d}{dx}{E}^{(k_1, k_2, \ldots, k_r)}_{n+1}(x;a,b,c)&=&(n+1)(r\ln c){E}^{(k_1, k_2, \ldots, k_r)}_{n}(x;a,b,c)\label{eqnmulti}
\end{eqnarray}
\end{thm}

The characterization of Appell polynomials is supposed to be used to establish the addition formula for the generalized multi poly-Euler polynomials using (\ref{eqnmulti}). But the constant $r\ln c$ that appeared in the expression of (\ref{eqnmulti}) disqualifies ${E}^{(k_1, k_2, \ldots, k_r)}_{n}(x;a,b,c)$ to be an Appell polynomial and to, consequently, satisfy any of the conditions of the said characterization. However, we can derive the addition formula using the same method in deriving the addition formula for ${E}^{(k_1, k_2, \ldots, k_r)}_{n}(x;a,b)$ in \cite{Jolany3}. The following theorem contains the addition formula for ${E}^{(k_1, k_2, \ldots, k_r)}_{n}(x;a,b,c)$.

\bigskip
\begin{thm}\label{multithm2}
The generalized poly-Euler polynomials satisfy the following addition formula
\begin{equation*}
{E}^{(k_1, k_2, \ldots, k_r)}_{n}(x+y;a,b,c)=\sum_{i=0}^n\binom{n}{i}(r\ln c)^{n-i}{E}^{(k_1, k_2, \ldots, k_r)}_{i}(x;a,b,c)y^{n-i}.
\end{equation*}

\begin{proof}
\begin{eqnarray*}
\sum_{n=0}^{\infty}{E}^{(k_1, k_2,\ldots, k_r)}_n(x+y;a,b,c)\frac{t^n}{n!}&=&\frac{2{\rm Li}_{(k_1, k_2,\ldots, k_r)}(1-(ab)^{-t})}{(a^{-t}+b^{t})^r}c^{(x+y)rt}\\
&=&\frac{2{\rm Li}_{(k_1, k_2,\ldots, k_r)}(1-(ab)^{-t})}{(a^{-t}+b^{t})^r}c^{xrt}c^{yrt}\\
&=&\left(\sum_{n=0}^{\infty}{E}^{(k_1, k_2,\ldots, k_r)}_n(x;a,b,c)\frac{t^n}{n!}\right)\left(\sum_{n=0}^{\infty}(yr\ln c)^n\frac{t^n}{n!}\right)\\
&=&\sum_{n=0}^{\infty}\left(\sum_{i=0}^n\binom{n}{i}(yr\ln c)^{n-i}{E}^{(k_1, k_2,\ldots, k_r)}_i(x;a,b,c)\right)\frac{t^n}{n!}.
\end{eqnarray*}
Comparing the coefficients of $\frac{t^n}{n!}$ yields the desired result.
\end{proof}
\end{thm}

\smallskip
The next theorem contains an explicit formula for ${E}^{(k_1, k_2, \ldots, k_r)}_{n}(x;a,b,c)$.

\bigskip
\begin{thm}\label{thm6}
The generalized multi poly-Euler polynomials have the following explicit formula
{\footnotesize
\begin{equation}\label{explicit3}
{E}^{(k_1,k_2,\ldots, k_r)}_n(x;a,b,c)=\sum_{i=0}^n\sum_{ 0\le m_1\le m_2\le\ldots\le m_r \atop c_1+c_2+\ldots=r}\sum_{j=0}^{m_r}\frac{2(rx\ln c-j\ln ab)^{n-i}r!(-1)^{j+s}(s\ln ab+r\ln a)^i\binom{m_r}{j}\binom{n}{i}}{(c_1!c_2!\ldots)(m_1^{k_1} m_2^{k_2}\ldots m_r^{k_r})},
\end{equation}
}
where $s=c_1+2c_2+\ldots$.

\begin{proof}
From Definition \ref{multigenpolyeuler}, we have
$${\rm Li}_{(k_1,k_2,\ldots, k_r)}(1-(ab)^{-t})c^{rxt}=\sum_{ 0\le m_1\le m_2\le\ldots\le m_r }\frac{(1-(ab)^{-t})^{m_r}}{m_1^{k_1} m_2^{k_2}\ldots m_r^{k_r}}e^{rxt\ln c}\qquad\qquad\qquad\qquad\qquad\qquad$$
\begin{eqnarray*}
&=&\sum_{ 0\le m_1\le m_2\le\ldots\le m_r }\frac{1}{m_1^{k_1} m_2^{k_2}\ldots m_r^{k_r}}\sum_{j=0}^{m_r}(-1)^j\binom{m_r}{j}\sum_{n=0}^{\infty}(rx\ln c-j\ln ab)^n\frac{t^n}{n!}\\
&=&\sum_{n=0}^{\infty}\left(\sum_{ 0\le m_1\le m_2\le\ldots\le m_r }\sum_{j=0}^{m_r}\frac{(-1)^{j}(rx\ln c-j\ln ab)^n\binom{m_r}{j}}{m_1^{k_1} m_2^{k_2}\ldots m_r^{k_r}}\right)\frac{t^n}{n!}.
\end{eqnarray*}
On the other hand, 
\begin{eqnarray*}
\left(\frac{1}{a^{-t}+b^{t}}\right)^r&=&a^{rt}\left(\frac{1}{1+(ab)^t}\right)^r=a^{rt}\left(\sum_{ n\ge0 }(-1)^n(ab)^{nt}\right)^r\\
&=&\sum_{c_1+c_2+\ldots=r}\frac{r!(-1)^{c_1+2c_2+\ldots}}{c_1!c_2!\ldots}e^{t[r\ln a+(c_1+2c_2+\ldots)\ln ab}\\
&=&\sum_{c_1+c_2+\ldots=r}\frac{r!(-1)^{c_1+2c_2+\ldots}}{c_1!c_2!\ldots}\sum_{n=0}^{\infty}(r\ln a+(c_1+2c_2+\ldots)\ln ab)^n\frac{t^n}{n!}\\
&=&\sum_{n=0}^{\infty}\left(\sum_{c_1+c_2+\ldots=r}\frac{r!(-1)^{c_1+2c_2+\ldots}(r\ln a+(c_1+2c_2+\ldots)\ln ab)^n}{c_1!c_2!\ldots}\right)\frac{t^n}{n!}.
\end{eqnarray*}
Hence,
$$\frac{2{\rm Li}_{(k_1,k_2,\ldots, k_r)}(1-(ab)^{-t})}{(a^{-t}+b^{t})^r}c^{rxt}=2{\rm Li}_{(k_1,k_2,\ldots, k_r)}(1-(ab)^{-t})e^{rxt\ln c}a^{rt}\left(\frac{1}{1+(ab)^{t}}\right)^r\qquad\qquad\qquad\qquad\qquad\qquad$$
\begin{eqnarray*}
&=&2\sum_{n=0}^{\infty}\sum_{i=0}^n\left(\sum_{ 0\le m_1\le m_2\le\ldots\le m_r }\sum_{j=0}^{m_r}\frac{(-1)^{j}(rx\ln c-j\ln ab)^{n-i}\binom{m_r}{j}}{m_1^{k_1} m_2^{k_2}\ldots m_r^{k_r}}\right)\frac{t^{n-i}}{(n-i)!}\times\\
&&\;\;\;\;\times\left(\sum_{c_1+c_2+\ldots=r}\frac{r!(-1)^{c_1+2c_2+\ldots}(r\ln a+(c_1+2c_2+\ldots)\ln ab)^i}{c_1!c_2!\ldots}\right)\frac{t^i}{i!}\\
&=&2\sum_{n=0}^{\infty}\sum_{i=0}^n\sum_{ 0\le m_1\le m_2\le\ldots\le m_r \atop c_1+c_2+\ldots=r}\sum_{j=0}^{m_r}\frac{H(r,i,j,n,a,b)}{(c_1!c_2!\ldots)(m_1^{k_1} m_2^{k_2}\ldots m_r^{k_r})}\frac{t^{n}}{n!}
\end{eqnarray*}
where 
$$H(r,i,j,n,a,b)=(rx\ln c-j\ln ab)^{n-i}r!(-1)^{j+c_1+2c_2+\ldots}(r\ln a+(c_1+2c_2+\ldots)\ln ab)^i\binom{m_r}{j}\binom{n}{i}.$$
By comparing the coefficient of $t^n/n!$, we obtain the desired explicit formula.
\end{proof}
\end{thm}

%%%%%%%%%%%%%%%%%%%%%%%%%%%%%%%

\begin{defn}\label{defn5}
For $m, n\ge0$, we define
\begin{equation}\label{symgen}
\mathcal{D}^{(m)}_n(x,y;a,b,c)=\sum_{k_1+k_2+\ldots +k_r=m}\binom{m}{k_1,k_2,\ldots k_r}\frac{{E}^{(-k_1,-k_2,\ldots -k_{r-1})}_n(x;a,b,c)}{(\ln a+\ln b)^n}\left(\frac{(r-1)y\ln c+\ln a}{\ln a +\ln b}\right)^{k_r}.
\end{equation}
\end{defn}

The following theorem contains the double generating function for $D^{(m)}_n(x,y;a,b,c)$.

\bigskip
\begin{thm}\label{thmm4}
For $n,m\ge0$, we have
\begin{equation}\label{eqnnnew6}
\sum_{n=0}^{\infty}\sum_{m=0}^{\infty}\mathcal{D}^{(m)}_n(x,y;a,b,c)\frac{t^n}{n!}\frac{u^m}{m!}=\frac{2e^{\left(\frac{(r-1)y\ln c+\ln a}{\ln a +\ln b}\right)u}e^{(r-1)\left(\frac{(r-1)x\ln c+\ln a}{\ln a +\ln b}\right)t}e^{\binom{r}{2}u+(r-1)t}(1-e^{-t})^{r-1}}{(1+e^{t})^{r-1}\prod_{i=1}^{r-1}(e^t+e^{iu}-e^{t+iu})}.
\end{equation}
\begin{proof}
$$\sum_{n=0}^{\infty}\sum_{m=0}^{\infty}\mathcal{D}^{(m)}_n(x,y;a,b,c)\frac{t^n}{n!}\frac{u^m}{m!}\qquad\qquad\qquad\qquad\qquad\qquad\qquad\qquad\qquad\qquad\qquad\qquad\qquad\qquad\qquad\qquad$$
\begin{eqnarray*}
&=&\sum_{n=0}^{\infty}\sum_{m=0}^{\infty}\sum_{k_1+k_2+\ldots +k_r=m}\frac{{E}^{(-k_1,-k_2,\ldots -k_{r-1})}_n(x;a,b,c)}{(\ln a+\ln b)^n}\left(\frac{(r-1)y\ln c+\ln a}{\ln a +\ln b}\right)^{k_r}\frac{t^n}{n!}\frac{u^m}{k_1!k_2!\ldots k_r!}\\
&=&\sum_{n=0}^{\infty}\sum_{k_1+k_2+\ldots +k_r\ge0}\frac{{E}^{(-k_1,-k_2,\ldots -k_{r-1})}_n(x;a,b,c)}{(\ln a+\ln b)^n}\left(\frac{(r-1)y\ln c+\ln a}{\ln a +\ln b}\right)^{k_r}\frac{t^n}{n!}\frac{u^{k_1+k_2+\ldots +k_r}}{k_1!k_2!\ldots k_r!}\\
&=&\sum_{n=0}^{\infty}\sum_{k_1+k_2+\ldots +k_{r-1}\ge0}\frac{{E}^{(-k_1,-k_2,\ldots -k_{r-1})}_n(x;a,b,c)}{(\ln a+\ln b)^n}\sum_{k_r\ge0}\left(\frac{(r-1)y\ln c+\ln a}{\ln a +\ln b}\right)^{k_r}\frac{u^{k_r}}{k_r!}\frac{t^n}{n!}\frac{u^{k_1+k_2+\ldots +k_{r-1}}}{k_1!k_2!\ldots k_{r-1}!}\\
&=&e^{\left(\frac{(r-1)y\ln c+\ln a}{\ln a +\ln b}\right)u}\sum_{n=0}^{\infty}\sum_{k_1+k_2+\ldots +k_{r-1}\ge0}\frac{{E}^{(-k_1,-k_2,\ldots -k_{r-1})}_n(x;a,b,c)}{(\ln a+\ln b)^n}\frac{t^n}{n!}\frac{u^{k_1+k_2+\ldots +k_{r-1}}}{k_1!k_2!\ldots k_{r-1}!}
\end{eqnarray*}
Using Theorem \ref{multithm1}, we obtain
$$\sum_{n=0}^{\infty}\sum_{m=0}^{\infty}\mathcal{D}^{(m)}_n(x,y;a,b,c)\frac{t^n}{n!}\frac{u^m}{m!}\qquad\qquad\qquad\qquad\qquad\qquad\qquad\qquad\qquad\qquad\qquad\qquad\qquad\qquad\qquad\qquad$$
\begin{eqnarray*}
&=&e^{\left(\frac{(r-1)y\ln c+\ln a}{\ln a +\ln b}\right)u}\sum_{k_1+k_2+\ldots +k_{r-1}\ge0}\sum_{n=0}^{\infty}{E}^{(-k_1,-k_2,\ldots -k_{r-1})}_n\left(\frac{(r-1)x\ln c+\ln a}{\ln a +\ln b}\right)\frac{t^n}{n!}\frac{u^{k_1+k_2+\ldots +k_{r-1}}}{k_1!k_2!\ldots k_{r-1}!}\\
&=&e^{\left(\frac{(r-1)y\ln c+\ln a}{\ln a +\ln b}\right)u}e^{(r-1)\left(\frac{(r-1)x\ln c+\ln a}{\ln a +\ln b}\right)t}\sum_{k_1+k_2+\ldots +k_{r-1}\ge0}\frac{2{\rm Li}_{(-k_1, -k_2,\ldots, -k_{r-1})}(1-e^{-t})}{(1+e^{t})^{r-1}}\frac{u^{k_1+k_2+\ldots +k_{r-1}}}{k_1!k_2!\ldots k_{r-1}!}\\
&=&\frac{2e^{\left(\frac{(r-1)y\ln c+\ln a}{\ln a +\ln b}\right)u}e^{(r-1)\left(\frac{(r-1)x\ln c+\ln a}{\ln a +\ln b}\right)t}}{(1+e^{t})^{r-1}}\sum_{0<m_1<m_2<\ldots <m_{r-1}}(1-e^{-t})^{m_{r-1}}\mathcal{S}(u,m_1,m_2,\ldots, m_{r-1})
\end{eqnarray*}
where
$$\mathcal{S}(u,m_1,m_2,\ldots, m_{r-1})=\sum_{k_1+k_2+\ldots +k_{r-1}\ge0}\frac{(um_1)^{k_1}(um_2)^{k_2}\ldots (um_{r-1})^{k_{r-1}}}{k_1!k_2!\ldots k_{r-1}!}\qquad\qquad\qquad\qquad\qquad\qquad\qquad\qquad\qquad$$
\begin{eqnarray*}
&=&\sum_{\widehat{m}\ge0}\frac{1}{\widehat{m}!}\sum_{k_1+k_2+\ldots +k_{r-1}=\widehat{m}}\binom{\widehat{m}}{k_1, k_2, \ldots k_{r-1}}(um_1)^{k_1}(um_2)^{k_2}\ldots (um_{r-1})^{k_{r-1}}\\
&=&\sum_{\widehat{m}\ge0}\frac{um_1+um_2+\ldots +um_{r-1})^{\widehat{m}}}{\widehat{m}!}\\
&=&e^{u(m_1+m_2+\ldots +m_{r-1})}.
\end{eqnarray*}
Thus,
$$\sum_{n=0}^{\infty}\sum_{m=0}^{\infty}\mathcal{D}^{(m)}_n(x,y;a,b,c)\frac{t^n}{n!}\frac{u^m}{m!}\qquad\qquad\qquad\qquad\qquad\qquad\qquad\qquad\qquad\qquad\qquad\qquad\qquad\qquad\qquad\qquad$$
\begin{eqnarray*}
&=&\frac{2e^{\left(\frac{(r-1)y\ln c+\ln a}{\ln a +\ln b}\right)u}e^{(r-1)\left(\frac{(r-1)x\ln c+\ln a}{\ln a +\ln b}\right)t}}{(1+e^{t})^{r-1}}\sum_{0<m_1<m_2<\ldots <m_{r-1}}(1-e^{-t})^{m_{r-1}}e^{u(m_1+m_2+\ldots +m_{r-1})}\\
&=&\frac{2e^{\left(\frac{(r-1)y\ln c+\ln a}{\ln a +\ln b}\right)u}e^{(r-1)\left(\frac{(r-1)x\ln c+\ln a}{\ln a +\ln b}\right)t}}{(1+e^{t})^{r-1}}\frac{e^u(1-e^{-t})}{1-e^u(1-e^{-t})}\frac{e^{2u}(1-e^{-t})}{1-e^{2u}(1-e^{-t})}\ldots\frac{e^{(r-1)u}(1-e^{-t})}{1-e^{(r-1)u}(1-e^{-t})}\\
&=&\frac{2e^{\left(\frac{(r-1)y\ln c+\ln a}{\ln a +\ln b}\right)u}e^{(r-1)\left(\frac{(r-1)x\ln c+\ln a}{\ln a +\ln b}\right)t}e^{\binom{r}{2}u}(1-e^{-t})^{r-1}}{(1+e^{t})^{r-1}\prod_{i=1}^{r-1}(1-e^{iu}(1-e^{-t}))}\\
&=&\frac{2e^{\left(\frac{(r-1)y\ln c+\ln a}{\ln a +\ln b}\right)u}e^{(r-1)\left(\frac{(r-1)x\ln c+\ln a}{\ln a +\ln b}\right)t}e^{\binom{r}{2}u+(r-1)t}(1-e^{-t})^{r-1}}{(1+e^{t})^{r-1}\prod_{i=1}^{r-1}(e^t+e^{iu}-e^{t+iu})}.
\end{eqnarray*}
\end{proof}
\end{thm}

\smallskip
Note that equation (\ref{eqnnew6}) can easily be deduced from equation (\ref{eqnnnew6}) by taking $r=1$. It is then interesting to establish an explicit formula for $\mathcal{D}^{(m)}_n(x,y;a,b,c)$ parallel to Theorem \ref{thm5}. To do this, let us consider first the following expression from the right-hand side of equation (\ref{eqnnnew6}). That is,
\begin{eqnarray*}
\frac{1}{(1+e^{t})^{r-1}\prod_{i=1}^{r-1}(e^t+e^{iu}-e^{t+iu})}&=&\left(\sum_{n=0}^{\infty}(-e)^{nt}\right)^{r-1}\prod_{i=1}^{r-1}\frac{1}{(1- (e^t-1)(e^{iu}-1))}\\
&=&\left(\sum_{n=0}^{\infty}(-e)^{nt}\right)^{r-1}\prod_{i=1}^{r-1}\sum_{j=0}^{\infty}(e^t-1)^j(e^{iu}-1)^j\\
&=&\left(\sum_{n=0}^{\infty}(-e)^{nt}\right)^{r-1}\prod_{i=1}^{r-1}\sum_{c_i\ge0}(e^t-1)^{c_i}(e^{iu}-1)^{c_i}\\
\left(\sum_{n=0}^{\infty}(-e)^{nt}\right)^{r-1}&=&\sum_{q=0}^{\infty}\sum_{k_{r-1}=0}^n\sum_{k_{r-2}=0}^{n-k_{r-1}}\ldots\sum_{k_{1}=0}^{n-k_{r-1}-\ldots -k_{2}}(-1)^qe^{qt}\\
&=&\sum_{q=0}^{\infty}(-1)^q\frac{\prod_{j=0}^{q-2}(q+1+j)}{(q-1)!}e^{qt}
\end{eqnarray*}
\begin{eqnarray*}
\prod_{i=1}^{r-1}\sum_{c_i\ge0}(e^t-1)^{c_i}(e^{iu}-1)^{c_i}&=&\sum_{j=0}^{\infty}\sum_{c_1+c_2+\ldots +c_{r-1}=j}(e^t-1)^{j}\prod_{i=1}^{r-1}(e^{iu}-1)^{c_i}\\
&=&\sum_{j=0}^{\infty}\sum_{c_1+c_2+\ldots +c_{r-1}=j}j!\sum_{n=0}^{\infty}\braced{n}{j}\frac{t^n}{n!}\prod_{i=1}^{r-1}c_i!\sum_{n=0}^{\infty}\braced{n}{c_i}\frac{u^n}{n!}
\end{eqnarray*}
\begin{eqnarray*}
\prod_{i=1}^{r-1}c_i!\sum_{n=0}^{\infty}\braced{n}{c_i}\frac{u^n}{n!}&=&\sum_{m=0}^{\infty}\sum_{d_{r-2}=0}^{m}\sum_{d_{r-3}=0}^{m-d_{r-2}}\ldots \sum_{d_{1}=0}^{m-d_{r-2}-\ldots - d_2}c_1!\braced{m-d_{r-2}-\ldots - d_{1}}{c_1}\times\\
&&\;\;\;\times\prod_{i=1}^{r-2}\binom{m-d_{r-2}-\ldots - d_{i+1}}{d_i}\braced{d_i}{c_{i+1}}c_{i+1}!(i+1)^{d_i}\frac{u^m}{m!}
\end{eqnarray*}
$$e^{(r-1)\left(\frac{(r-1)x\ln c+\ln a}{\ln a +\ln b}\right)t}(1-e^{-t})^{r-1}\left(\sum_{n=0}^{\infty}(-e)^{nt}\right)^{r-1}\left(j!\sum_{n=0}^{\infty}\braced{n}{j}\frac{t^n}{n!}\right)
\qquad\qquad\qquad\qquad\qquad\qquad\qquad$$
\begin{eqnarray*}
&=&e^{(r-1)\left(\frac{(r-1)x\ln c+\ln a}{\ln a +\ln b}\right)t}\left(\sum_{k=0}^{r-1}(-1)^k\binom{r-1}{k}e^{-kt}\right)\left(\sum_{q=0}^{\infty}(-1)^q\frac{\prod_{j=0}^{q-2}(q+1+j)}{(q-1)!}e^{qt}\right)\left(j!\sum_{n=0}^{\infty}\braced{n}{j}\frac{t^n}{n!}\right)\\
&=&\left(\sum_{q=0}^{\infty}\sum_{k=0}^{r-1}e^{\left(\frac{(r-1)^2x\ln c+(q-k)\ln b+(q-k+r-1)\ln a}{\ln a+\ln b}\right)}\frac{(-1)^{k+q}\binom{r-1}{k}\prod_{j=0}^{q-2}(q+1+j)}{(q-1)!}\right)\left(j!\sum_{n=0}^{\infty}\braced{n}{j}\frac{t^n}{n!}\right)\\
&=&\left(\sum_{n=0}^{\infty}\sum_{q=0}^{\infty}\sum_{k=0}^{r-1}\left(\frac{(r-1)^2x\ln c+(q-k)\ln b+(q-k+r-1)\ln a}{\ln a+\ln b}\right)^{n}\right.\\
&&\;\;\;\;\;\left.\frac{(-1)^{k+q}\binom{r-1}{k}\prod_{j=0}^{q-2}(q+1+j)}{(q-1)!}\frac{t^n}{n!}\right)\left(j!\sum_{n=0}^{\infty}\braced{n}{j}\frac{t^n}{n!}\right)\\
&=&\sum_{n=0}^{\infty}\left(\sum_{p=0}^n\sum_{q=0}^{\infty}\sum_{k=0}^{r-1}\binom{n}{p}\left(\frac{(r-1)^2x\ln c+(q-k)\ln b+(q-k+r-1)\ln a}{\ln a+\ln b}\right)^{n-p}\right.\\
&&\;\;\;\;\;\left.\frac{(-1)^{k+q}\binom{r-1}{k}\prod_{j=0}^{q-2}(q+1+j)}{(q-1)!}j!\braced{p}{j}\right)\frac{t^n}{n!}
\end{eqnarray*}
$$e^{\left(\frac{(r-1)y\ln c+\ln a}{\ln a +\ln b}\right)u}e^{\binom{r}{2}u}\left(\prod_{i=1}^{r-1}c_i!\sum_{n=0}^{\infty}\braced{n}{c_i}\frac{u^n}{n!}\right)
\qquad\qquad\qquad\qquad\qquad\qquad\qquad\qquad\qquad\qquad\qquad$$
\begin{eqnarray*}
&=&\left(\sum_{m=0}^{\infty}\left(\frac{(r-1)y\ln c+\binom{r}{2}\ln b+\left\{\binom{r}{2}+1\right\}\ln a}{\ln a +\ln b}\right)^m\frac{u^m}{m!}\right)\left(\sum_{m=0}^{\infty}\sum_{d_{r-2}=0}^{m}\sum_{d_{r-3}=0}^{m-d_{r-2}}\ldots\right.\\
&&\;\;\;\left. \sum_{d_{1}=0}^{m-d_{r-2}-\ldots - d_2}c_1!\braced{m-d_{r-2}-\ldots - d_{1}}{c_1}\prod_{i=1}^{r-2}\binom{m-d_{r-2}-\ldots - d_{i+1}}{d_i}\braced{d_i}{c_{i+1}}c_{i+1}!(i+1)^{d_i}\frac{u^m}{m!}\right)\\
&=&\sum_{m=0}^{\infty}\left\{\sum_{l=0}^m\left(\frac{(r-1)y\ln c+\binom{r}{2}\ln b+\left\{\binom{r}{2}+1\right\}\ln a}{\ln a +\ln b}\right)^{m-l}\sum_{d_{r-2}=0}^{l}\sum_{d_{r-3}=0}^{l-d_{r-2}}\ldots\right.\\
&&\;\;\;\left. \sum_{d_{1}=0}^{l-d_{r-2}-\ldots - d_2}c_1!\braced{l-d_{r-2}-\ldots - d_{1}}{c_1}\prod_{i=1}^{r-2}\binom{l-d_{r-2}-\ldots - d_{i+1}}{d_i}\braced{d_i}{c_{i+1}}c_{i+1}!(i+1)^{d_i}\right\}\frac{u^m}{m!}
\end{eqnarray*}
$$\sum_{n=0}^{\infty}\sum_{m=0}^{\infty}\mathcal{D}^{(m)}_n(x,y;a,b,c)\frac{t^n}{n!}\frac{u^m}{m!}\qquad\qquad\qquad\qquad\qquad\qquad\qquad\qquad\qquad\qquad\qquad\qquad\qquad\qquad\qquad\qquad$$
\begin{eqnarray*}
&=&\sum_{n=0}^{\infty}\sum_{m=0}^{\infty}\left\{\sum_{j=0}^{\infty}\sum_{c_1+c_2+\ldots +c_{r-1}=j}\sum_{l=0}^m\left(\frac{(r-1)y\ln c+\binom{r}{2}\ln b+\left\{\binom{r}{2}+1\right\}\ln a}{\ln a +\ln b}\right)^{m-l}\right.\\
&&\;\;\;\left. \sum_{d_{r-2}=0}^{l}\sum_{d_{r-3}=0}^{l-d_{r-2}}\ldots\sum_{d_{1}=0}^{l-d_{r-2}-\ldots - d_2}c_1!\braced{l-d_{r-2}-\ldots - d_{1}}{c_1}\prod_{i=1}^{r-2}\binom{l-d_{r-2}-\ldots - d_{i+1}}{d_i}\times\right.\\
&&\;\;\left. \times\braced{d_i}{c_{i+1}}c_{i+1}!(i+1)^{d_i}\sum_{p=0}^n\sum_{q=0}^{\infty}\sum_{k=0}^{r-1}\binom{n}{p}\left(\frac{(r-1)^2x\ln c+(q-k)\ln b+(q-k+r-1)\ln a}{\ln a+\ln b}\right)^{n-p}\right.\\
&&\;\;\;\;\;\left.\frac{(-1)^{k+q}\binom{r-1}{k}\prod_{j=0}^{q-2}(q+1+j)}{(q-1)!}j!\braced{p}{j}\right\}\frac{t^n}{n!}\frac{u^m}{m!}
\end{eqnarray*}
Comparing coefficients, we obtain the following theorem.

\bigskip
\begin{thm}\label{thm555}
For $n,m\ge0$, we have
$$\mathcal{D}^{(m)}_n(x,y;a,b,c)=\sum_{j=0}^{\infty}\sum_{c_1+c_2+\ldots +c_{r-1}=j}\sum_{l=0}^m\left(\frac{(r-1)y\ln c+\binom{r}{2}\ln b+\left\{\binom{r}{2}+1\right\}\ln a}{\ln a +\ln b}\right)^{m-l}\times$$
\begin{eqnarray*}
&&\;\times\sum_{d_{r-2}=0}^{l}\sum_{d_{r-3}=0}^{l-d_{r-2}}\ldots\sum_{d_{1}=0}^{l-d_{r-2}-\ldots - d_2}c_1!\braced{l-d_{r-2}-\ldots - d_{1}}{c_1}\prod_{i=1}^{r-2}\binom{l-d_{r-2}-\ldots - d_{i+1}}{d_i}\times\\
&&\;\times\braced{d_i}{c_{i+1}}c_{i+1}!(i+1)^{d_i}\sum_{p=0}^n\sum_{q=0}^{\infty}\sum_{k=0}^{r-1}\binom{n}{p}\left(\frac{(r-1)^2x\ln c+(q-k)\ln b+(q-k+r-1)\ln a}{\ln a+\ln b}\right)^{n-p}\times\\
&&\;\times\frac{(-1)^{k+q}\binom{r-1}{k}\prod_{j=0}^{q-2}(q+1+j)}{(q-1)!}j!\braced{p}{j}.
\end{eqnarray*}
\end{thm}

\bigskip
\begin{flushleft}
{\bf Hassan Jolany}\\
Universit\'e des Sciences et Technologies de Lille\\
UFR de Math\'ematiques\\
Laboratoire Paul Painlev\'e\\
CNRS-UMR 8524 59655 Villeneuve d'Ascq Cedex/France\\
e-mail: hassan.jolany@math.univ-lille1.fr
\end{flushleft}

\bigskip
\begin{flushleft}
{\bf Roberto B. Corcino}\\
Mathematics and ICT Department\\
Cebu Normal University\\
Osmena Blvd., Cebu City\\
Philippines 6000\\
e-mail: rcorcino@yahoo.com
\end{flushleft}

\end{document}